*Reducing time-dependent multifactor Black and Scholes equation with knock-out features to equivalent time-constant coefficient equation, and applications*

*Tark Bouhennache, PhD.*

*Abstract*

We consider the multifactor Black and Scholes equation with time-dependent coefficients, and a knock-out feature contingent on the underlying asset values reaching a limit (reflected by a Dirichlet condition on the boundary). We prove that this equation, which has important applications in finance and insurance, can be reduced to an equivalent time-constant coefficient equation, with coefficients defined as averages of the original ones. Equivalent results are also valid for general second order parabolic equations, with applications in other fields in the natural sciences. The result established in this article has not been documented so far in the presence of boundary conditions. The proof is provided in a general framework, as it invokes techniques from the Functional Analysis theory, namely the Hille-Yosida approximation technique. A main ingredient in the proof is establishing the $e^{-A_1-A_2} = e^{-A_1}e^{-A_2}$ identity, proved here in a general setting for unbounded, commuting, monotone and maximal operators $A_1, A_2$. This would in particular allow generalization of the result to more general classes of evolution problems with time dependent generators.

*1. Introduction*

We consider the following multifactor Black and Scholes equation (B&S equation)

$$\frac{\partial U}{\partial \tau} + \sum_{i=1}^{n}(r-m)y_i \frac{\partial U}{\partial y_i} + \frac{1}{2}\sum_{i,j=1}^{n}\rho_{i,j}\sigma_i\sigma_j y_i y_j \frac{\partial^2 U}{\partial y_i \partial y_j} - rU - dU = 0, \qquad \forall \tau < T, y \in \Omega', \tag{1}$$

with a Dirichlet boundary condition, and a final condition at $\tau = T$. We assume that the coefficients $r, m, d, \sigma_i$ and $\rho_{i,j}, i,j \leq n$, in the differential operator above are piecewise continuous with respect to the time variable $\tau$ and independent of the space variable $y = (y_1, \ldots, y_n) \in \Omega'$, where $\Omega'$ is an open subset of $\mathbb{R}_+^n$ with smooth boundary. We assume the matrix $\left(\rho_{i,j}(\tau)\sigma_i(\tau)\sigma_j(\tau)\right)_{ij}$ to be symmetric and definite-positive, uniformly for $\tau \in [0,T]$, so that there exists a constant $c > 0$ such that:

$$\sum_{i,j \leq n}\rho_{i,j}(\tau)\sigma_i(\tau)\sigma_j(\tau)\xi_i\xi_j \geq c \sum_{i \leq n}\xi_i^2, \qquad \forall \xi \in \mathbb{R}^n \text{ and } \forall \tau \in [0,T]. \tag{2}$$

In this article we present a general proof that, given the final condition $U(T,.)$, one can compute $U(\tau_0,.)$, for any given $\tau_0 \in [0,T]$, using only the time-averages of the coefficients in (1) over the interval $[\tau_0, T]$. More precisely, let us consider $\overline{U}$ to be the solution of the following equivalent time-constant coefficients equation:

$$\frac{\partial \bar{U}}{\partial \tau} + \sum_{i=1}^{n}(\bar{r}-\bar{m})y_i \frac{\partial \bar{U}}{\partial y_i} + \frac{1}{2}\sum_{i,j=1}^{n} \bar{\rho}_{i,j} \bar{\sigma}_i \bar{\sigma}_j y_i y_j \frac{\partial^2 \bar{U}}{\partial y_i \partial y_j} - \bar{r}\,\bar{U} - \bar{d}\bar{U} = 0, \qquad \forall \tau < T, y \in \Omega', \qquad (3)$$

with the same Dirichlet conditions and final condition $\bar{U}(T,y) = U(T,y), \forall y \in \Omega'$, and where

$$\bar{r} = \frac{1}{T-\tau_0}\int_{\tau_0}^{T} r(s)ds, \quad \bar{m} = \frac{1}{T-\tau_0}\int_{\tau_0}^{T} m(s)ds, \quad \bar{\sigma}_i^2 = \frac{1}{T-\tau_0}\int_{\tau_0}^{T} \sigma_i(s)^2 ds, \qquad (4)$$

$$\bar{d} = \frac{1}{T-\tau_0}\int_{\tau_0}^{T} d(s)ds \text{ and } \bar{\rho}_{i,j} = \frac{1}{(T-\tau_0)\bar{\sigma}_i\bar{\sigma}_j}\int_{\tau_0}^{T} \sigma_i(s)\sigma_j(s)\rho_{i,j}(s)ds. \qquad (5)$$

We will prove that $\bar{U}(\tau_0,.) = U(\tau_0,.)$. This is the main result of this article, and so referred to throughout this document. This result has not been documented so far in the presence of boundary conditions in the general case. The Dirichlet boundary condition for both equations above is $U(\tau,y) = 0$, for all $y \in \partial\Omega'$, with $y \notin \partial\mathbb{R}_+^n$.

Reducing the one factor B&S equation to a constant coefficient equation has been studied in the literature in the case where $\Omega' = \mathbb{R}_+$, and we can refer to ] and [9] where simple change of variables are suggested for the proof. Instead, one can use more mathematical but straightforward proofs, and in the more general multifactor case (with $\Omega' = \mathbb{R}_+^n$). For the sake of completeness we will present two such proofs in the appendix. The first proof uses the Risk Neutral valuation approach that establishes the solution of (1) in terms of an expectation value with respect to a Risk Neutral measure. The second uses the Fourier transform, which requires $\Omega' = \mathbb{R}_+^n$. These proofs however do not carry over to the general case where $\Omega'$ is any open subset of $\mathbb{R}_+^n$ with smooth boundary, and it is the contribution of this article to present a proof in this more general case. The latter invokes techniques from Functional Analysis, namely the Hille-Yosida theory.

For convenience the presentation, the main result will be re-stated and proved in terms of second order parabolic equations with initial condition. In fact, using the well-known change of variables: $x_i = \text{Ln}(y_i), i = 1, \ldots, n$, and $t = (T-\tau)$, and denoting $u(t,x) = U(\tau,y)$, transforms Equation (1) as follows:

$$\frac{\partial u}{\partial t} + \mathcal{A}(t) = 0, \text{ with } \mathcal{A}(t) = -\frac{1}{2}\sum_{i,j \leq n}\frac{\partial}{\partial x_i}\left(\rho_{i,j}\sigma_i\sigma_j\frac{\partial u}{\partial x_j}\right) + \sum_{i \leq n}\left(\frac{\sigma_i^2}{2} - (r-m)\right)\frac{\partial u}{\partial x_i} + (r+d)u, \qquad (6)$$

with the Dirichlet condition on the boundary. The coefficients of operator $\mathcal{A}(t)$ are independent on the space variable $x = (x_1, \ldots, x_n)$. Note that the equation above also shows that this paper's main result has applications in fields such as physics and engineering, in addition to the financial applications as we emphasize below. It is also worth mentioning that the general proof presented, based on Functional Analysis arguments, could also be used to establish the result for more general classes of evolution problems.

For the more interested reader more background on the finance interpretation of equation (1) is provided in the Appendix. In the case where $d = 0$, the B&S equation (1) models the price of European type financial derivative instruments, such as put and call options, on a basket of assets. The assets are assumed to be correlated, with their vector values denoted by $S(\tau) = (S_1(\tau), \ldots, S_n(\tau))$. The Dirichlet condition on the boundary of $\Omega'$ reflects a knock-out feature by which the derivative instrument loses value and expires whenever $S(\tau)$ hits the boundary of $\Omega'$. More generally, where $d$ is not necessarily zero, the B&S equation (1) arises in the Fair Market Valuation and risk management of investment guarantees embedded in variable annuity products. Reducing equation (1) to the time-constant coefficient case simplifies solving the equation, and in particular allows for the solution to be expressed analytically in certain situations, for example when using the famous B&S formulas. The parameter $d$ reflects policy decrements, more specifically mortality and lapse decrements. Although in practice lapse rates are "dynamic", in the sense that $d$ also depends on the space variable $y$, we are assuming in this article that all parameters in equation (1), including $d$, are independent of $y$. The latter condition is necessary in establishing the main result.

The general proof of the main result of this article is organized in three steps as follows:

**Step 1:** As mentioned we re-state the B&S equation in terms of a second order parabolic equation: $\frac{\partial u}{\partial t} + \mathcal{A}(t)u = 0$, on $[0,T]X\Omega$, for a certain domain $\Omega$, with Dirichlet condition on the boundary and an initial condition at $t = 0$.

**Step 2:** we prove the result when the coefficients of $\mathcal{A}(t)$ are piecewise constant. As we will see this is a direct consequence of the $e^{-A_1-A_2} = e^{-A_1}e^{-A_2}$ identity that we will prove in a *general setting* for commuting, monotone and maximal operators $A_1, A_2$ with the property that $A_1 + A_2$ is monotone and maximal. The proof uses techniques from the Hille-Yosida theory, see [3] for more details on that theory.

**Step 3**: we prove and use a convergence result to show that the result established in step 2 for operators $\mathcal{A}(t)$ with piecewise constant coefficients carries over to the more general case where the coefficients of the operator are piecewise continuous with respect to the time variable $t$.

This article is organized as follows. In section 2, we redefine the result in terms of second order parabolic equations $\frac{\partial u}{\partial t} + \mathcal{A}(t)u = 0$, and state the main theorem. In section 3 we prove the key identity $e^{-A_1-A_2} = e^{-A_1}e^{-A_2}$. Section 4 proceeds with the proof of the main theorem with the key ingredients being the $e^{-A_1-A_2} = e^{-A_1}e^{-A_2}$ identity, along with a convergence result that will be proved. In the Appendix we provide, for the more interested reader and for completeness of the presentation, more details on the financial background of the B&S equation (1), and in particular on the knock-out feature of financial derivative instruments. We also present two direct and short overviews of proofs, with less attention to mathematical rigor, of the main result in the case where $\Omega' = \mathbb{R}_+^n$, which are more mathematical and general proofs than those in ] and [9].

**Acknowledgements**: My sincere gratitude goes to Professor James Ralston from UCLA for his support and encouragements, and for kindly providing feedback on the manuscript. I also wish to thank Mourad Sini from the Austrian Academy of Science who kindly accepted to review and provide helpful feedback.

## *2. Statement of the main result*

For ease of presentation we will assume that $\Omega'$ is *strictly* included in $\mathbb{R}_+^n$, in which case we will assume in the remainder of this article, and without loss of generality, that the Dirichlet condition holds on the whole boundary $\partial \Omega'$. Without loss of generality we also set $\tau_0 = 0$.

As shown above, the known change of variables $x_i = \text{Ln}(y_i), i = 1, \ldots, n,$ and $t = T - \tau$, transform the B&S equation (1) with Dirichlet condition into the second order parabolic equation (6) with a Dirichlet condition on the boundary of the domain $\Omega \stackrel{\text{def}}{=} \{x = (x_1, \ldots, x_n) \in \mathbb{R}^n; (e^{x_1}, \ldots, e^{x_n}) \in \Omega'\}$. We will be stating and proving the main result in terms of second order parabolic problems of the following form:

$$\begin{cases} \frac{\partial u}{\partial t} + \mathcal{A}(t)u = 0, & \forall t < T, x \in \Omega, \\ u(t,x) = 0, & \forall t < T, x \in \partial\Omega, \\ u(0,x) = g(x), & \forall x \in \Omega. \end{cases} \quad (7)$$

For each $t \geq 0$, the operator $\mathcal{A}(t)$ as defined in equation (6) is a second order elliptic operator of the form:

$$\mathcal{A}(t)u = -\sum_{i,j\leq n}\frac{\partial}{\partial x_i}\left(a_{i,j}(t)\frac{\partial u}{\partial x_j}\right)+\sum_{i\leq n}b_i(t)\frac{\partial u}{\partial x_i}+q(t)u. \tag{8}$$

In our case we have $a_{i,j}(t) = \rho_{i,j}(T-t)\sigma_i(T-t)\sigma_j(T-t)$, with $b_i$ and $q$ defined in a similar way. From (2) this operator is then elliptic uniformly with respect to $t \in [0,T]$, in the sense that there exists a constant $c > 0$ such that:

$$\sum_{i,j\leq n}a_{i,j}(t)\xi_i\xi_j \geq c\sum_{i\leq n}\xi_i^2, \quad \forall\,\xi \in \mathbb{R}^n \text{ and } \forall t \in [0,T]. \tag{9}$$

The matrix $\left(a_{i,j}(t)\right)_{i,j\leq n}$ is also symmetric (and definite positive). In addition the coefficients $a_{i,j}(t), b_i(t)$ and $q(t)$ of the operator $\mathcal{A}(t)$ are piecewise continuous functions of the variable $t$, so uniformly bounded on the interval $[0,T]$, and are independent of the space variable $x$.

Before we proceed, let us introduce some notations and assumptions. As a reminder, $L^2(\Omega)$ is the Hilbert space of measurable and real-valued functions $u$ defined on $\Omega \subset \mathbb{R}^n$ (not to confuse the generic notation $u$ with that of the solution to the parabolic equation above) equipped with the following scalar product:

$$<u,v> = \int_\Omega u(x)v(x)\,dx, \text{ and norm } \|u\|_{L^2(\Omega)} \stackrel{\text{def}}{=} \left(\int_\Omega u(x)^2\,dx\right)^{1/2}.$$

Also, $H_0^1(\Omega)$ is the Hilbert space of real-valued functions $u$ such that $u$ and $\frac{\partial u}{\partial x_i} \in L^2(\Omega), \forall i \leq n$, with $u$ satisfying the Dirichlet condition on the boundary of $\Omega$. Here the derivatives $\frac{\partial u}{\partial x_i}$ are defined in the weak sense. $H_0^1(\Omega)$ is equipped with the following norm: $\|u\|_{H_0^1(\Omega)} \stackrel{\text{def}}{=} \left(\int_\Omega \left(u(x)^2 + \sum_{i\leq n}\frac{\partial u(x)^2}{\partial x_i}\right)dx\right)^{1/2}$. For any given Hilbert space, or a Banach space, $B$ the space $L^2(0,T;B)$ (resp. $L^\infty(0,T;B)$) is that of measurable functions $u$ defined on the interval $[0,T]$ with values in $B$ and equipped with the following norm:

$$\|u\|_{L^2(0,T;B)} \stackrel{\text{def}}{=} \left(\int_0^T \|u(t)\|_B^2\,dt\right)^{\frac{1}{2}}, \quad \left(\text{resp. } \|u\|_{L^\infty(0,T;B)} = \max_{t\in[0,T]}\|u(t)\|_B\right).$$

Some of the spaces $B$ that we will be considering are $L^2(\Omega)$, $H_0^1(\Omega)$ or $H^{-1}(\Omega)$. The space $H^{-1}(\Omega)$ is the dual space of $H_0^1(\Omega)$. For more information on Hilbert spaces and Sobolev spaces we refer to [3] and [4].

Let us denote by $a_t(u,v)$, for any given $t \geq 0$, the following bilinear form $(u,v) \to a_t(u,v) \in \mathbb{R}$:

$$a_t(u,v) = \int_\Omega\left(\sum_{i,j\leq n}a_{i,j}(t)\frac{\partial u(x)}{\partial x_i}\frac{\partial v(x)}{\partial x_j}+\sum_{i\leq n}b_i(t)\frac{\partial u(x)}{\partial x_i}v(x)+q(t)u(x)v(x)\right)dx, \tag{10}$$

$\forall u,v \in H_0^1(\Omega)$. It is bounded on $H_0^1(\Omega)\times H_0^1(\Omega)$, uniformly for $t \in [0,T]$, in the sense that there exists a constant $c_0 > 0$ such that $|a_t(u,u)| \leq c_0\|u\|_{H_0^1(\Omega)}^2$, $\forall u,v \in H_0^1(\Omega)$ and $t \in [0,T]$. Thanks to (9) the latter bilinear form is also coercive, or equivalently bounded from below, in the sense that there exist constants $c_1, c_2 > 0$, such that:

$$|a_t(u,u)| + c_1\|u\|_{L^2(\Omega)}^2 \geq c_2\|u\|_{H_0^1(\Omega)}^2, \quad \forall u \in H_0^1(\Omega), \quad \forall t \in [0,T]. \tag{11}$$

The operator $\mathcal{A}(t)$ is linked to the form $a_t(.,.)$ in the sense that for each $u \in H_0^1(\Omega)$ such that $\mathcal{A}(t)u \in L^2(\Omega)$, and using integrations by part, we have

$$a_t(u,v) = <\mathcal{A}(t)u, v>, \quad \forall v \in H_0^1(\Omega).$$

A function $u(.,.)$ of the variables $t$ and $x$ can be interpreted as a function of $t$ taking values in a space of functions of $x \in \Omega$, such that $[u(t)](x) = u(t,x)$. Let us now give a reminder on the notion of weak solutions of the following inhomogeneous second order parabolic problem (as opposed to the homogeneous equation (7) where the term $f$ is set equal to zero):

$$\begin{cases} \dfrac{\partial u}{\partial t} + \mathcal{A}(t)u = f, & \forall t < T, x \in \Omega, \\ u(t,x) = 0, & \forall t < T, x \in \partial\Omega, \\ u(0,x) = g(x), & \forall x \in \Omega. \end{cases} \quad (12)$$

We refer to [4] for more details, in particular on existence and uniqueness results for parabolic equations.

**_Definition_**: Given $g \in L^2(\Omega)$ and $f \in L^2(0,T; H^{-1}(\Omega))$, a function $u \in L^2(0,T; H_0^1(\Omega))$, with the derivative $\dfrac{du}{dt} \in L^2(0,T; H^{-1}(\Omega))$, is a weak solution of the second order parabolic problem above if and only if

$$<\frac{du(t,.)}{dt}, v>_{H^{-1}, H_0^1} + a_t(u(t,.), v) = <f(t,.), v>_{H^{-1}, H_0^1}, \quad \forall v \in H_0^1(\Omega), \text{ with } u(0) = g. \quad (13)$$

As indicated in [4] the weak solution is necessarily continuous: $u \in C(0,T; L^2(\Omega))$, so the initial condition $u(0,.) = g(.)$ is well defined.

In situations that arise in practice the initial condition $g$ is not necessarily square integrable, which is the case for example for a Put option where $g(x) = MAX(K - e^x, 0)$ is bounded but not square integrable. To allow for such situations we will work with "larger" functional spaces and consider the function $g$ to be possibly an element of $L^{2,s}(\Omega)$, for a certain $s \in \mathbb{R}$, where $L^{2,s}(\Omega)$ is by definition the Hilbert space of measured real-valued functions equipped with the following scalar product:

$$<u,v>_{L^{2,s}(\Omega)} \stackrel{\text{def}}{=} \int_\Omega u(x)v(x)(1+|x|^2)^s \, dx, \text{ with } \|u\|^2_{L^{2,s}(\Omega)} \stackrel{\text{def}}{=} \int_\Omega u(x)^2(1+|x|^2)^s \, dx.$$

Here $|x|^2 \stackrel{\text{def}}{=} x_1^2 + \cdots + x_n^2$. We also define $H_0^{1,s}(\Omega)$ to be the Hilbert space of real-valued functions satisfying the Dirichlet condition on the boundary of $\Omega$ and equipped with the following norm

$$\|u\|^2_{H_0^{1,s}(\Omega)} \stackrel{\text{def}}{=} \int_\Omega \left(u(x)^2 + \sum_{i \le n} \frac{\partial u(x)^2}{\partial x_i}\right)(1+|x|^2)^s \, dx.$$

We denote by $H^{-1,-s}(\Omega)$ the dual space of $H_0^{1,s}(\Omega)$, which is the space of bounded linear forms defined on $H_0^{1,s}(\Omega)$. The space $L^{2,-s}(\Omega)$ is in particular a subset of $H^{-1,-s}(\Omega)$. Similar to (13), a weak solution $u$ to the parabolic equation (12) with initial condition $g \in L^{2,s}(\Omega)$, and with $f \in L^2(0,T; H^{-1,s}(\Omega))$, is a function $u \in L^2(0,T; H_0^{1,s}(\Omega))$, with $\frac{du}{dt} \in L^2(0,T; H^{-1,s}(\Omega))$, such that the equality in (13) holds for all $\forall v \in H_0^{1,-s}(\Omega)$. We have the following existence and uniqueness theorem.

**_Theorem 1_** (*existence of weak solutions*): For every $f \in L^2(0,T; H^{-1,s}(\Omega))$ and $g \in L^{2,s}(\Omega)$, there exists a unique weak solution of the problem (12). In addition the solution satisfies the following energy estimates:

$$\max_{0\le t\le T} \|u(t)\|^2_{L^{2,s}(\Omega)} + \|u\|^2_{L^2(0,T;H^{1,s}_0(\Omega))} + \left\|\frac{du}{dt}\right\|^2_{L^2(0,T;H^{-1,s}(\Omega))}$$
$$\le C\left(\|f\|^2_{L^2(0,T;H^{-1,s}(\Omega))} + \|g\|^2_{L^{2,s}(\Omega)}\right), \tag{14}$$

*for a certain constant $C$. These results also hold if the coefficients of the operator (8) are smooth functions of the space variable $x$ and the estimate (9) holds uniformly for $t$ and $x$.*

**Proof**. We first show that if the theorem holds for $s = 0$, then it holds for any $s \in \mathbb{R}$. Let us denote by $\mathcal{L}_s$ the multiplying operator, such that $\mathcal{L}_s u = (1 + |x^2|)^{-s/2} u$. The operator $\mathcal{L}_s$ is obviously a linear, bounded and invertible operator from $L^2(\Omega)$ onto $L^{2,s}(\Omega)$. It is also bounded from $H^1_0(\Omega)$ onto $H^{1,s}_0(\Omega)$, and from $H^{-1}(\Omega)$ onto $H^{-1,s}(\Omega)$.

From operator $\mathcal{A}(t)$ we can derive the operator $\mathcal{B}(t) = \mathcal{L}_s^{-1}\mathcal{A}(t)\mathcal{L}_s$, which we can easily show to be elliptic, with bounded coefficient and with second order coefficients satisfying (9) uniformly in $t$ and $x$. Let $\hat{g} = \mathcal{L}_s^{-1}g \in L^2(\Omega)$ and $\hat{f} = \mathcal{L}_s^{-1}f \in L^2(0,T;H^{-1}(\Omega))$, then based on the assumption that the theorem holds for $s = 0$ there exists a solution $\hat{u}(t) \in L^2(0,T;H^1_0(\Omega))$ of the parabolic equation $\frac{\partial \hat{u}}{\partial t} + \mathcal{B}(t)\hat{u} = \hat{f}$, with initial condition $\hat{u}(0) = \hat{g}$. The solution $u$ for $s \ne 0$ is then given by $\mathcal{L}_s\hat{u}(t)$.

The proof of the theorem when $s = 0$ is detailed in [4] where it is based on constructing a Galerkin approximation to the solution of the form

$$u_m(t,x) = \sum_{k=1}^m d_m^k(t) w_k(x),$$

where $w_k, k = 1, \dots, +\infty$, is an orthogonal basis in $H^1_0(\Omega)$. The only difference here is that we are not assuming $\Omega$ to be bounded since this condition is not necessary for the existence of such an orthogonal basis. This completes the proof. **Q.E.D.**

Let us now define, for any given $T > 0$, the following second order elliptic operator with time-constant coefficients:

$$\frac{1}{T}\int_0^T \mathcal{A}(s)ds \stackrel{\text{def}}{=} \bar{\mathcal{A}}_T \quad \text{where} \quad \bar{\mathcal{A}}_T u \stackrel{\text{def}}{=} \sum_{i,j\le n}\frac{\partial}{\partial x_i}\left(\bar{a}_{i,j}\frac{\partial u}{\partial x_j}\right) + \sum_{i\le n}\bar{b}_i\frac{\partial u}{\partial x_j} + \bar{q}u, \quad \forall u,$$

derived by taking the following average values of the coefficients of operator $\mathcal{A}(t)$ as for (4):

$$\bar{a}_{i,j} \stackrel{\text{def}}{=} \frac{1}{T}\int_0^T a_{i,j}(s)ds, \quad \bar{b}_i \stackrel{\text{def}}{=} \frac{1}{T}\int_0^T b_i(s)ds, \quad \text{and} \quad \bar{q} \stackrel{\text{def}}{=} \frac{1}{T}\int_0^T q(s)ds. \tag{15}$$

Let us denote by $\bar{u}_T$ the weak solution of the following second order parabolic problem *with constant coefficients*:

$$\begin{cases} \dfrac{\partial \bar{u}_T}{\partial t} + \bar{\mathcal{A}}_T\bar{u}_T = 0, & \forall t < T, x \in \Omega, \\ \bar{u}_T(t,x) = 0, & \forall t < T, x \in \partial\Omega, \\ \bar{u}_T(0,x) = g(x), & \forall x \in \Omega. \end{cases} \tag{16}$$

The following is the main result of this article:

**_Theorem 2 (main)_**: Let $g \in L^{2,s}(\Omega)$, for a given $s \in \mathbb{R}$. Denote by $\boldsymbol{u}$ the weak solution of (7), and by $\bar{\boldsymbol{u}}_T$ the weak solution of (16). Then under the conditions above $\boldsymbol{u}(T) = \bar{\boldsymbol{u}}_T(T)$, for any given $T > 0$.

**Highlight of the proof:** The proof is based on the main lemma that is stated and proved in the next section. It then proceeds in three steps as outlined below:

1/ In **step 1** we use the main lemma to prove the theorem in the case where the operator-valued function $t \to \mathcal{A}(t)$ is piecewise constant and $g \in H_0^1(\Omega) \cap H^2(\Omega)$.

2/ In **step 2** we prove a convergence result and use it to demonstrate the theorem when the coefficients of the operator-valued function $t \to \mathcal{A}(t)$ are piecewise continuous functions.

3/ In **step 3** we use density results to complete the proof of the theorem for any $g \in L^{2,s}(\Omega)$.

Note that in the case where the coefficients of the differential operator $\mathcal{A}(t)$ are not dependent on the time variable, existence and uniqueness also result from the general Hille-Yosida theory, see [3] for example.

A required condition under this theory is that $A(t) + c_1$, for a certain constant $c_1$, is a monotone and maximal operator. In general, we denote by $A(t)$ the restriction of $\mathcal{A}(t)$ on the domain:

$$D(A(t)) \stackrel{\text{def}}{=} \{u \in H_0^1(\Omega); \exists f \in L^2(\Omega), \text{ such that } a_t(u,v) = <f,v>, \forall v \in H_0^1(\Omega)\}. \quad (17)$$

Using integration by parts we have: $<A(t)u, v> = a_t(u,v)$, $\forall u \in D(A(t)), \forall v \in H_0^1(\Omega)$. Since the second order elliptic operator $\mathcal{A}(t)$ has coefficients that do not depend on $x$, we can show that $D(A(t)) = D(A) \stackrel{\text{def}}{=} H_0^1(\Omega) \cap H^2(\Omega)$ (which holds in general if the coefficients of $A(t)$ depend smoothly on $x$, and given that $\Omega$ has a smooth boundary). In particular, the latter domain does not depend on $t$. Given (11) we have that

$$<A(t)u, u> + c_1 \parallel u \parallel_{L^2(\Omega)}^2 \geq 0 \quad \forall u \in D(A(t)),$$

which means that the operator $A(t) + c_1$ is monotone. On the other hand, given that the form $a_t(.,.)$ is coercive, or equivalently bounded from below, the associated operator $A(t)$ is such that for any $z > c_1$ the operator $A(t) + z$ defined from $D(A(t))$ onto $L^2(\Omega)$ is invertible (see [4], chapter 6. Or also [6], chapter 6 - similar argument is used for an equation stemming from physical science in [1]), which implies that the operator $A(t) + c_1$ is maximal as well. Without loss of generality we will assume in the proof that $c_1 = 0$.

When $A \stackrel{\text{def}}{=} A(t)$ is not dependent on $t$, the Hille-Yosida theorem states that for any $g \in D(A)$, there exists a continuous function with values in $D(A)$ such that $\frac{du}{dt} + Au = 0$, with $u(0) = g$. This solution is then necessarily equal to the weak solution established in *Theorem* 1. It is usually expressed in the exponential form: $u(t) = e^{-tA}g$, *thus giving sense to exponential of unbounded operator*.

Note that a direct consequence of the differential operator $\mathcal{A}(t)$ not depending on the space variable is the commutation property $A(t_1)A(t_2) = A(t_2)A(t_1)$, for any $t_1, t_2$, which is essential in our proof of the main result. If we assume that $t \to A'(t)$ is a certain continuous function with values in the space of *bounded linear operators* defined on $L^2(\Omega)$ (equipped with the corresponding norm), then the differential equation $\frac{du}{dt} + A'(t)u = 0$, with initial value $u(0) = g$, has a unique solution which is infinitely differentiable with respect to $t$. If in addition the commutation property $A'(t_1)A'(t_2) = A'(t_2)A'(t_1)$ holds, for any $t_1, t_2$, then the solution is given by $u(t) =$ $e^{-\int_0^t A'(s)ds}g$, with the integral defined in the Riemann sense. *Theorem* 2 can then be seen as a generalization for unbounded operators of this exponential expression of the solution.

## 3. Main lemma: the identity $e^{-A_1-A_2} = e^{-A_1}e^{-A_2}$ for unbounded commuting operators

For any $t_1, t_2$, the operators $A_1 = A(t_1)$ and $A_2 = A(t_2)$, as defined above, satisfy the following properties:

**1/** Commutation property: $A_1 A_2 = A_2 A_1$.

**2/** Operators $A_1$, $A_2$ and $A_1 + A_2$ are monotone and maximal (a property that can be relaxed by assuming it holds after possibly adding a constant).

**3/** Operators $A_1$ and $A_2$ have the same domain that we will denote by $D(A)$, so $D(A) = D(A_1) = D(A_2)$, and have equivalent graph norms. Recall that for any operator $A$, the graph norm is defined as follows:

$$\| u \|_{D(A)} \stackrel{\text{def}}{=} (\| u \|_H^2 + \| Au \|_H^2)^{\frac{1}{2}}. \tag{18}$$

As we mentioned before, in the case of operators $A(t)$ derived from elliptic operators not depending on the space variable, the graph norm above is equivalent to that $H^2(\Omega)$. The following result is *key* in the proof of the main theorem, which pertains to operators derived from elliptic second order operators, but the proof of the lemma below will be *provided in a general framework* where $A_1$ and $A_2$ are operators defined on any Hilbert space $H$ equipped with a scalar product $<.,.>_H$ and satisfying the three properties above.

**Lemma 3 (Main)**: If operators $A_1$ and $A_2$ satisfy properties 1, 2 and 3 above, then $e^{-A_1-A_2} = e^{-A_1}e^{-A_2}$.

The exponential of unbounded operators can be defined in different ways, and in this article we define it as mentioned based on the Hille-Yosida theory. More specifically, let $t \to u(t)$ be the solution of the following evolution problem:

$$\begin{cases} \dfrac{du}{dt} + Au = 0 & \text{on } [0, +\infty), \\ u(0) = u_0 \in D(A), \end{cases} \tag{19}$$

under the condition that the unbounded and linear operator $A: D(A) \subset H \to H$ is monotone and maximal, in the sense that

$$< Av, v > \geq 0, \forall v \in D(A), \text{ and } I + \lambda A: D(A) \subset H \to H \text{ is bijective for a certain given } \lambda \geq 0.$$

The exponential of operator $A$ is such that $u(t) = e^{-tA}u_0$ is the solution of (19). Before we proceed with the proof we will first provide a brief reminder of some useful results from the Hille-Yosida theory. See [3] for more information.

**Reminder of some results from the Hille-Yosida theory:**

**Theorem 4** *(Existence and uniqueness)*: For each $u_0 \in D(A)$, problem (19) has a unique solution $u(t) = e^{-tA}u_0$ in $C^1([0, +\infty); H) \cap C^0([0, +\infty); D(A))$, where the space $D(A)$ is equipped with the graph norm (18). In addition, the following inequalities hold:

$$\| e^{-tA}u_0 \|_H \leq \| u_0 \|_H \quad \text{and} \quad \left\| \frac{de^{-tA}u_0}{dt} \right\|_H = \| Au \|_H \leq \| Au_0 \|_H \quad \forall t \geq 0. \tag{20}$$

The result above holds for bounded operators. For an unbound operator $A$ the proof as presented in [3] is based viewing $A$ as a limit of Yosida approximation $A_\lambda$, as $\lambda \to 0$, where:

$$A_\lambda = \frac{1}{\lambda}(I - J_\lambda), \quad \text{with } J_\lambda = (I + \lambda A)^{-1}. \tag{21}$$

The operator $A_\lambda$ satisfies the following properties:

- It is a bounded, monotone and maximal operator,
- $\| A_\lambda v \|_H \leq \| Av \|_H$ and $A_\lambda v \to A(v)$, as $\lambda \to 0, \forall v \in D(A)$.

As a reminder the proof of the theorem above as presented in [3] proceeds by introducing, for each $\lambda > 0$, the following evolution problem:

$$\begin{cases} \dfrac{du_\lambda(t)}{dt} + A_\lambda u_\lambda(t) = 0 & \text{on } [0, +\infty), \\ u_\lambda(0) = u_0, \end{cases}$$

for which exists a solution, since $A_\lambda$ is a linear and *bounded* operator. The following inequalities are then established:

$$\| u_\lambda(t) - u_\mu(t) \|_H^2 \leq 4(\lambda + \mu)t \ \| Au_0 \|_H^2, \quad \forall u_0 \in D(A),$$

and

$$\left\| \frac{du_\lambda(t)}{dt} - \frac{du_\mu(t)}{dt} \right\|_H^2 \leq 4(\lambda + \mu)t \ \| A^2 u_0 \|_H^2, \quad \forall u_0 \in D(A^2).$$

which show uniform convergence of $u_\lambda$ and $\frac{du_\lambda}{dt}$, as $\lambda \to 0$, to a limit $u \in C^1([0, +\infty); H)$ when $u_0 \in D(A^2)$. As presented in [3], using the density of $D(A^2)$ in $D(A)$, with $D(A^2) = \{w \in D(A) | Aw \in D(A)\}$, it is then proved that results in the theorem above are also valid when $u_0 \in D(A)$.

Note that a consequence of $A$ being monotone and maximal is that $D(A)$ is dense in $H$, and $A$ is a closed operator. Since $D(A^2) = (I + A)^{-1} D(A)$ and $(I + A)^{-1}$ is bounded in $H$, $D(A^2)$ is also dense in $H$.

***Note on the Regularity of the solution $u(t)$:*** One can show that the operator $A$ is also monotone and maximal when restricted on the Hilbert space $H' = D(A)$, equipped with the following scalar product:

$$<u, v>_{D(A)} = <u, v>_H + <Au, Av>_H, \tag{22}$$

with domain $D(A^2)$ as defined above. One can prove that $H'$ is a closed space based on the fact that $A$ is a closed operator, and is thus a Hilbert space. As a consequence, *Theorem* 4 can also be applied on $H'$ when $u_0 \in D(A^2)$ to show that $u(t) \in C^1([0, +\infty); D(A)) \cap C^0([0, +\infty); D(A^2))$.

**Proof of *Lemma* 3 (Main):** Since operators $A_1$ and $A_2$ have equivalent norms of the graph, there exist constants $a_1, a_2 > 0$ such that $\| A_1 v \|_H \leq a_2 \| A_2 v \|_H$ and $\| A_2 v \|_H \leq a_1 \| A_1 v \|_H, \forall v \in D(A) \stackrel{\text{def}}{=} D(A_1) = D(A_2)$. The proof of the Lemma is based on the use of the Yosida operator approximations $A_{i,\lambda}$ of the operators $A_i, i = 1, 2$, as defined in (21). We will show that for every $u_0 \in D(A)$:

1/ $e^{-A_{1,\lambda}} e^{-A_{2,\lambda}} u_0 \to e^{-A_1} e^{-A_2} u_0$, as $\lambda \to 0$, and

2/ $e^{-A_{1,\lambda} - A_{2,\lambda}} u_0 \to e^{-A_1 - A_2} u_0$, as $\lambda \to 0$.

The Lemma then results from the fact that $e^{-A_{1,\lambda}} e^{-A_{2,\lambda}} = e^{-A_{1,\lambda} - A_{2,\lambda}}$ since the operators $A_{1,\lambda}$ and $A_{2,\lambda}$ are bounded and commuting (which is consequence of $A_1$ and $A_2$ commuting). Note that the identity $e^{-A_{1,\lambda}} e^{-A_{2,\lambda}} = e^{-A_{1,\lambda} - A_{2,\lambda}}$ can be directly proved for bounded operators using the power series expansion definition of the exponential of an operator, similar to that for the exponential of numbers. To complete the proof of the lemma we proceed in 4 steps as follows:

**Step 1:** Let us prove the convergence $e^{-A_{1,\lambda}} e^{-A_{2,\lambda}} u_0 \to e^{-A_1} e^{-A_2} u_0$, as $\lambda \to 0$, for any given $u_0 \in D(A)$. We have

$$e^{-A_{1,\lambda}} e^{-A_{2,\lambda}} u_0 - e^{-A_1} e^{-A_2} u_0 = e^{-A_{1,\lambda}} (e^{-A_{2,\lambda}} u_0 - e^{-A_2} u_0) + (e^{-A_{1,\lambda}} - e^{-A_1}) e^{-A_2} u_0.$$

Note that according to Hille-Yosida theorem $e^{-A_2} u_0$ also belongs to $D(A)$, which shows the convergence to zero of the second term in the right hand side in the equality above as $\lambda \to 0$. The convergence to zero of the first term is consequence of the following inequality: $\| e^{-A_{1,\lambda}} (e^{-A_{2,\lambda}} u_0 - e^{-A_2} u_0) \|_H \leq \| e^{-A_{2,\lambda}} u_0 - e^{-A_2} u_0 \|_H$, since we $A_{1,\lambda}$ is monotone.

**Step 2:** In this step we prove the convergence of $e^{-t(A_{1,\lambda} + A_{2,\lambda})} u_0$, as $\lambda \to 0$, to a certain function $\bar{u}(t)$. The convergence as we will see is uniform with respect to $t$ in any bounded interval. Let us initially assume that $u_0 \in D(A^2)$. It is important to note that $A_{1,\lambda} + A_{2,\lambda}$ is not the Yosida approximation of $A_1 + A_2$, but a sum of two operators that are Yosida approximations, so the convergence is not a direct consequence from the Hille-Yosida theory. Nevertheless, we show that *similar techniques* apply to prove the convergence. By definition, $u_\lambda(t) = e^{-t(A_{1,\lambda} + A_{2,\lambda})} u_0$ is the unique solution to the following evolution problem:

$$\begin{cases} \dfrac{du_\lambda}{dt} + (A_{1,\lambda} + A_{2,\lambda}) u_\lambda = 0, & \text{on } [0, +\infty), \\ u_\lambda(0) = u_0. \end{cases} \quad (23)$$

Following similar techniques as in [3], we have

$$\frac{1}{2} \frac{d}{dt} \| u_\lambda - u_\mu \|_H^2 + \sum_{i=1}^{2} < A_{i,\lambda} u_\lambda - A_{i,\mu} u_\mu, u_\lambda - u_\mu >_H = 0, \quad (24)$$

and for each of $i = 1,2$, the following inequality holds:

$$< A_{i,\lambda} u_\lambda - A_{i,\mu} u_\mu, u_\lambda - u_\mu >_H \geq < A_{i,\lambda} u_\lambda - A_{i,\mu} u_\mu, \lambda A_{i,\lambda} u_\lambda - \mu A_{i,\mu} u_\mu >_H.$$

We also have, since $< t(A_{1,\lambda} + A_{2,\lambda}) v, v >_H \geq 0, \forall v \in H$, that $t \to \| e^{t(A_{1,\lambda} + A_{2,\lambda})} v \|_H$ is a decreasing function of $t$.

Since the operators $A_{1,\lambda}$ and $A_{2,\lambda}$ commute, using the power series definition of the exponential, we have

$$A_{i,\lambda} u_\lambda(t) = A_{i,\lambda} e^{t(A_{1,\lambda} + A_{2,\lambda})} u_0 = e^{t(A_{1,\lambda} + A_{2,\lambda})} A_{i,\lambda} u_0.$$

It then follows that $\| A_{i,\lambda} u_\lambda(t) \|_H \leq \| A_{i,\lambda} u_0 \|_H \leq \| A_i u_0 \|_H$. From (24) we then get

$$\frac{1}{2}\frac{d}{dt}\|u_\lambda - u_\mu\|_H^2 \leq 4(\lambda + \mu)(\|A_1 u_0\|_H^2 + \|A_2 u_0\|_H^2),$$

which, by integrating, yields the following inequality

$$\|u_\lambda(t) - u_\mu(t)\|_H \leq 2\sqrt{2(\lambda + \mu)t}\,(\|A_1 u_0\|_H^2 + \|A_2 u_0\|_H^2)^{\frac{1}{2}}, \tag{25}$$

showing that $u_\lambda$ is uniformly convergent to a function $\bar{u}$ on every bounded interval $[0, T']$, $\forall\, T' > 0$.

This result holds for every $u_0 \in D(A)$. Since we assumed in this step that $u_0 \in D(A^2)$ (here $D(A^2) \stackrel{\text{def}}{=} D(A_1^2) = D(A_2^2)$), we can in addition show that function $v_\lambda = \frac{du_\lambda}{dt}$ is also uniformly convergent, which proves that $\bar{u} \in C^1[0, +\infty; H]$, $\forall\, T' > 0$, and that $\frac{du_\lambda}{dt}$ converges to $\frac{d\bar{u}}{dt}$, as $\lambda \to 0$.

Note that $v_\lambda$ satisfies the first equality in equation (23) with initial condition $v_\lambda(0) = -(A_{1,\lambda} + A_{2,\lambda})u_0$. Similar to (25) we get the following inequality for a certain constant $c$:

$$\|v_\lambda(t) - v_\mu(t)\|_H \leq 4\sqrt{(\lambda + \mu)t}\,(\|A_1 v_\lambda(0)\|_H^2 + \|A_2 v_\mu(0)\|_H^2)^{\frac{1}{2}}$$
$$\leq c\sqrt{(\lambda + \mu)t}\,(\|A_1^2 u_0\|_H^2 + \|A_2^2 u_0\|_H^2)^{\frac{1}{2}},$$

since $u_0 \in D(A^2)$, and where we have used the fact that there exists a constant $c'$ such that

$$\|A_{i,\lambda} v_\lambda(t)\|_H \leq \|A_{i,\lambda} v_\lambda(0)\|_H \leq \|A_i v_\lambda(0)\|_H \leq c'(\|A_1^2 u_0\|_H^2 + \|A_2^2 u_0\|_H^2)^{\frac{1}{2}},$$

This then shows the uniform convergence of the derivative function $v_\lambda = \frac{du_\lambda}{dt}$.

**Step 3:** In this step we still assume $u_0 \in D(A^2)$, and prove under this condition that the limit function $\bar{u}(t)$ is equal to $e^{-t(A_1 + A_2)}u_0$, which proves the convergence $e^{-A_{1,\lambda} - A_{2,\lambda}}u_0 \to e^{-A_1 - A_2}u_0$, as $\lambda \to 0$, by setting $t = 1$. For that purpose we only need to show that the function $\bar{u}(t)$ satisfies the following equation:

$$\begin{cases} \dfrac{d\bar{u}}{dt} + (A_1 + A_2)\bar{u} = 0 & \text{on } [0, +\infty), \\ \bar{u}(0) = u_0. \end{cases} \tag{26}$$

Given that $\frac{du_\lambda}{dt} + (A_{1,\lambda} + A_{2,\lambda})u_\lambda = 0$, and taking the limit, as $\lambda \to 0$, we see that the first term $\frac{du_\lambda}{dt}$ in this equation converges to $\frac{d\bar{u}}{dt}$. For the second term, let us first note that the space $H' = D(A)$, equipped with the scalar product (22), is also closed and is then a Hilbert space. The restriction of operator $A = A_1 + A_2$ on $H'$, with domain $D(A^2)$, is also monotone and maximal. Since we assumed $u_0 \in D(A^2)$, the inequality (25) then also applies with the norm of $H' = D(A)$, which proves the uniform convergence of $u_\lambda$ in the graph norm $D(A)$.

Let us consider now the following decomposition:

$$(A_{1,\lambda} + A_{2,\lambda})u_\lambda - (A_1 + A_2)\bar{u} = \big((A_{1,\lambda} + A_{2,\lambda})u_\lambda - (A_1 + A_2)u_\lambda\big) + \big((A_1 + A_2)u_\lambda - (A_1 + A_2)\bar{u}\big), \tag{27}$$

which converges to zero for each $t$ as we show below. In fact, for the first term in the right hand side of this equality we have:

$$\big((A_{1,\lambda} + A_{2,\lambda})e^{-t(A_{1,\lambda} + A_{2,\lambda})}u_0 - (A_1 + A_2)e^{-t(A_{1,\lambda} + A_{2,\lambda})}u_0\big) = e^{-t(A_{1,\lambda} + A_{2,\lambda})}\big((A_{1,\lambda} + A_{2,\lambda}) - (A_1 + A_2)\big)u_0,$$

so

$$\left\|\left((A_{1,\lambda} + A_{2,\lambda})u_\lambda(t) - (A_1 + A_2)u_\lambda(t)\right)\right\|_H \leq \left\|(A_{1,\lambda} + A_{2,\lambda})u_0 - (A_1 + A_2)u_0\right\|_H,$$

which converges to zero. For the second term in the right hand side of equation (27), we note that $u_\lambda$ uniformly converges in the norm of $D(A)$, as we just established above, which proves that $\left((A_1 + A_2)u_\lambda - (A_1 + A_2)\bar{u}\right)$ converge to zero for the norm of $H$.

**Step 4**: The results in step 3 were established assuming that $u_0 \in D(A^2)$. Since $D(A^2)$ is dense in $D(A)$, due to the fact that $A$ is a closed operator as shown in [3], there exists a sequence $u_n \in D(A^2), n > 0$, such that $u_n \to u_0$, as $n \to +\infty$.

From the first inequality in (20) from the Hille-Yosida theorem, and since both $A_1$ and $A_1$ are monotone, we have that:

$$\| e^{-A_1}e^{-A_2}(u_n - u_0) \|_H \leq \| e^{-A_2}(u_n - u_0) \|_H \leq \| (u_n - u_0) \|_H,$$

which converges to zero as $n \to +\infty$.

On the other hand, the inequality $\| e^{-A_1-A_2}(u_n - u_0) \|_H \leq \| (u_n - u_0) \|_H$ holds since $A_1 + A_2$ is a monotone operator. The first term in the latter inequality then converges to zero, as $n \to +\infty$. Since $e^{-A_1}e^{-A_2}u_n = e^{-A_1-A_2}u_n$, for all $u_n \in D(A^2)$, we can then take the limit as $n \to +\infty$, to prove that $e^{-A_1}e^{-A_2}u_0 = e^{-A_1-A_2}u_0$, for all $u_0 \in D(A)$, which concludes the proof of the Lemma. **Q.E.D.**

In this article Lemma 3 will be applied to an operator $A = A_1 + A_2$ that is maximal since related to a second order elliptic operator. Establishing conditions on $A_1$ and $A_2$ for $A = A_1 + A_2$ to be maximal in the general case is a difficult problem, but one which would help generalizing the main result of this paper to more general classes of evolution problems. We have the following proposition:

**_Proposition_**: *Let $A_1$ and $A_2$ be two monotone and maximal operators on a Hilbert space $H$ such that $D(A_1) = D(A_2)$. The operator $A = A_1 + A_2$ is monotone and maximal under the assumption that there exist constants $a, b > 0$ such that the following inequality holds uniformly for $\alpha \in [0,1]$:*

$$\| A_1u \|_H + \| A_2u \|_H \leq b \| A_1u + \alpha A_2u \|_H + a \| u \|_H, \qquad \forall u \in D(A). \tag{28}$$

**Proof**: First, note that inequality above implies in particular that $A_1$ (respectively $A_2$) is $A_2$-bounded (respectively $A_1$-bounded). Also, note that it is straightforward to show that the following inequality holds for any integer $n > 0$:

$$\frac{1}{n} \| A_2u \|_H \leq \frac{b'}{n} \| A_1u + \alpha A_2u + u \|_H + \frac{a'}{n} \| u \|_H, \qquad \forall u \in D(A), \tag{29}$$

uniformly for $\alpha \in [0,1]$, for some constants $a', b' > 0$.

Using inequality (28) with $\alpha = 1$ and the fact that both $A_1$ and $A_2$ are closed, we show that $A = A_1 + A_2$ is also closed. To show that it is maximal, it suffices to show that $I + A$ is invertible. For that purpose, we proceed by induction showing that $I + A_1 + \frac{m}{n}A_2$ is invertible for all $m = 0, \ldots, n$, for a certain integer $n > 0$ selected to be large enough. For $m = 0$ the result is true since $I + A_1$ is invertible, given that $A_1$ is monotone and maximal.

We will use the following theorem from [6] (page 196): *If $T$ is a closed and invertible operator with a bounded inverse, and $B$ a $T$-bounded operator such as $\| Bu \|_H \leq c \| u \|_H + d \| Tu \|_H$ and $c \| T^{-1} \| + d < 1$, then $T + B$ is closed and invertible with a bounded inverse.*

Assume for a certain integer $m$, such as $0 < m < n$, that $I + A_1 + \frac{m}{n} A_2$ is invertible. Now with $n$ large enough we can apply the previous theorem to $T = I + A_1 + \frac{m}{n} A_2$ and $B = \frac{A_2}{n}$, using inequality (29), to show that operator $I + A_1 + \frac{m+1}{n} A_2$ is also invertible, thus completing the proof by induction. In applying the theorem above we have used the fact that $\| T^{-1} \| \leq 1$, which is a consequence of $A_1$ and $A_2$ being monotone operators.

Since $I + A_1 + \frac{m}{n} A_2$ is invertible for all $m \leq n$, we set $m = n$ to then show that $I + A_1 + A_2$ is invertible, which proves that $A_1 + A_2$ is maximal. **Q.E.D.**

## 4. Proof of the main theorem

We now return to the proof of the main theorem, which, as highlighted in section 3, proceeds in three steps as follows. As a reminder $A(t)$ is the restriction of the second order elliptic differential operator $\mathcal{A}(t)$, defined in (8), on the domain $D(A(t))$. The latter is defined in (17), which as mentioned is equal to $H_0^1(\Omega) \cap H^2(\Omega)$.

**Step 1**: we prove the result assuming that the operator-valued function $t \to A(t)$ is piecewise constant, which equivalently means that the coefficients of the operator $\mathcal{A}(t)$ are piecewise constant in the time variable $t$. We also assume that the initial condition $g$ belongs to $D(A(t)) = H_0^1(\Omega) \cap H^2(\Omega)$.

Let $0 = t_0 < t_1, \ldots, < t_{N-1} < t_N = T$, for a certain integer $N > 0$. Assume that $t \to A(t)$ takes one of $N$ operator-values: $A(t_i), \ i = 0, \ldots, N - 1$. For each $t \geq 0$, we assume that the operator $A(t)$ is defined as follows:

$$A(t) = \begin{cases} A(t_i), & \text{for } i \text{ such that } t \in [t_i, t_{i+1}[, \\ A(t_{N-1}), & \text{for } t \geq t_{N-1}. \end{cases}$$

Let $u$ be the solution of the evolution problem (7) where the operator-valued function $t \to A(t)$ is assumed to be piecewise constant. Applying the Hille-Yosida theorem on each separate interval $[t_i, t_{i+1}[$, with appropriate initial condition on each interval as we specify below, shows the existence of a solution $u(t) \in C(0, +\infty; D(A))$ such that $u \in C^1([t_i, t_{i+1}[; H)$, for each $i$.

For $i = 0$, the initial condition on the interval $[t_i, t_{i+1}[$ is set to be $u(0) = g$. For $i > 0$, the initial condition is selected to be the final limit value at $t_i$ of the solution solved on the preceding interval $[t_{i-1}, t_i[$. Note that the derivative $\frac{du}{dt}$ has a limit on the left at each $t_i$ that is not necessarily equal to $\frac{du(t_i)}{dt}$, so the function $\frac{du}{dt}$ could be discontinuous a $t = t_i$. The solution $u(T)$ at time $t = t_N = T$ is then given by:

$$u(T) = e^{-(T-t_{N-1})A(t_{N-1})} \ldots e^{-(t_2-t_1)A(t_1)} e^{-t_1 A(t_0)} g = \prod_{i=0}^{N-1} e^{-(t_{N-i}-t_{N-1-i})A(t_{N-1-i})} g \tag{30}$$

On the other hand, let us introduce the following operator:

$$\bar{A}_T = \frac{1}{T} \int_0^T A(t) dt \stackrel{\text{def}}{=} \sum_{i=0}^{N-1} \frac{t_{i+1} - t_i}{T} A(t_i),$$

which is monotone and maximal since it is associated with a coercive bilinear form as that in (10). Applying the Hille-Yosida theorem for the following evolution problem then shows the existence and uniqueness of the solution denoted by $\bar{u}_T$:

$$\begin{cases} \dfrac{d\bar{u}_T}{dt} + \bar{A}_T \bar{u}_T = 0 & \text{on } [0, +\infty), \\ \bar{u}_T(0) = g. \end{cases}$$

Let us use the induction argument to show that $u(T) = \bar{u}_T(T) \stackrel{\text{def}}{=} e^{-\int_0^T A(t)dt} g$. For $N = 2$ we have $u(T) = e^{-(t_2-t_1)A(t_1)} e^{-t_1 A(t_0)} g$ from (30). Applying Lemma 3 (main), knowing that operators $A(t_0)$ and $A(t_1)$ satisfy the three conditions of the lemma as announced at the beginning of section 4, the latter is then equal to

$$u(T) = e^{-(t_2-t_1)A(t_1) - t_1 A(t_0)} u_0 = e^{-T\left(\left(\frac{t_2-t_1}{T}\right)A(t_1) + \frac{t_1}{T}A(t_0)\right)} g = e^{-T\bar{A}_T} g = \bar{u}(T).$$

Assuming that the result is true for $N$, we follow the same argument to show that it is true for $N + 1$ too, which completes the proof. **Q.E.D.**

**Step 2**: we now prove the result when the coefficients of the operator $A(t)$ are piecewise continuous functions of $t$. We base the proof on the result from step 1 above, along with the convergence result from the following lemma.

Let $v$ be the solution to the parabolic equation (7). Let $\mathcal{A}_N(t), N \geq 0$, be a sequence of elliptic operators of the form

$$\mathcal{A}_N(t)u = \sum_{i,j \leq n} \frac{\partial}{\partial x_i}\left(a_{i,j}^N(t) \frac{\partial u}{\partial x_j}\right) + \sum_{i \leq n} b_i^N(t) \frac{\partial u}{\partial x_i} + q^N(t)u,$$

with coefficients uniformly converging to those of $\mathcal{A}(t)$ in the sense that, for each $i, j \leq n$, we have:

$$\max_{0 \leq t \leq T} \left(|a_{i,j}(t) - a_{i,j}^N(t)| + |b_i(t) - b_i^N(t)| + |q(t) - q^N(t)|\right) \to 0, \text{ as } N \to +\infty.$$

The matrix $\left(a_{i,j}^N(t)\right)_{ij}$ is assumed to be symmetric and uniformly definite-positive, in the sense that the inequality (9) is satisfied, uniformly for $t \in [0, T]$ in bounded intervals. Let $v_N$ be the solution to the parabolic equation (7) with the operator $\mathcal{A}(t)$ in the equation replaced by $\mathcal{A}_N(t)$. We have the following lemma.

<u>**Lemma 5**</u> *(convergence)*: Under the conditions above, assuming that $g \in L^2(\Omega)$, we have that

$$\max_{0 \leq t \leq T} \| v(t) - v_N(t) \|_{L^2(\Omega)} \to 0, \qquad \text{as } N \to +\infty.$$

**Proof**: The function $w_N = v - v_N$ is solution to the following equations

$$\begin{cases} \dfrac{\partial w_N}{\partial t} + \mathcal{A}(t)w_N = f_N, & \forall t < T, x \in \Omega, \\ w_N(t, x) = 0, & \forall t < T, x \in \partial\Omega, \\ w_N(0, x) = 0, & \forall x \in \Omega. \end{cases}$$

where $f_N = \left(\mathcal{A}(t) - \mathcal{A}_N(t)\right)v_N$. Since $g \in L^2(\Omega)$, then from *Theorem 1*, for a.e. $t \in [0, T]$, $w_N(t)$ belongs to $H_0^1(\Omega)$, which implies that $f_N(t) \in H^{-1}(\Omega)$, and we have

$$\| f_N \|_{L^2(0,T; H^{-1}(\Omega))} \leq c_N \| v_N \|_{L^2(0,T; H_0^1(\Omega))} \leq c_N c' \| g \|_{L^2(\Omega)}^2,$$

where there exists a constant $\bar{c} > 0$, such that

$$c_N \leq \bar{c} \max_{i,j \leq n} \left( \max_{0 \leq t \leq T} \left( |a_{i,j}(t) - a_{i,j}^N(t)| + |b_i(t) - b_i^N(t)| + |q(t) - q^N(t)| \right) \right),$$

which converges to zero, as $N \to +\infty$. From *Theorem* 1 we have, for a certain constant $c'' > 0$,

$$\max_{0 \leq t \leq T} \| w_N(t) \|_{L^2(\Omega)} \leq c'' \| f_N \|_{L^2(0,T; H^{-1}(\Omega))} \leq c_N c' c'' \| g \|_{L^2(\Omega)}^2,$$

which then converges to zero and proves the lemma. **Q.E.D.**

To complete the proof in step 2 we assume without loss of generality that $g \in H_0^1(\Omega) \cap H^2(\Omega)$, since the latter is dense in $L^2(\Omega)$. In fact, for any $h \in L^2(\Omega)$, there exists a sequence $g_m \in H_0^1(\Omega) \cap H^2(\Omega)$, $m \geq 0$, that converges to $h$ in $L^2(\Omega)$. If we prove that the result is true for each $g_m$, we will then take the limit, as $m \to +\infty$, to prove the result for $h \in L^2(\Omega)$. We used here the first inequality in (14) that shows that the solution of the parabolic problem depends "continuously" on the initial value.

For any given integer $N > 0$, let $t_k = \frac{kT}{N}$, $k = 0, \ldots, N$, such that $0 = t_0 < t_1 < \ldots < t_{N-1} < t_N = T$ is a uniform subdivision of the interval $[0, T]$. For each $i, j \leq n$, let us define the following piecewise constant functions:

$$a_{i,j}^N(t) = a_{i,j}(t_k), \quad b_i^N(t) = b_i(t_k) \text{ and } q^N(t) = q(t_k), \quad \forall t \in [t_k, t_{k+1}[, \forall k \leq N - 1.$$

Without loss of generality, we assume that all coefficients of operator $\mathcal{A}(t)$ are continuous on the left with respect to $t$, and that all points of possible discontinuity are included in the set $\{t_1, \ldots, t_{N-1}\}$. Let $\mathcal{A}^N(t)$ be the second order elliptic operator with the above coefficients, and let $u^N$ be the solution of the parabolic equation (7) with the operator $\mathcal{A}(t)$ in the equation replaced by $\mathcal{A}^N(t)$. Since the coefficients of the operator $\mathcal{A}(t)$ are piecewise continuous function, we have:

$$\max_{0 \leq t \leq T} \left( |a_{i,j}(t) - a_{i,j}^N(t)| + |b_i(t) - b_i^N(t)| + |q(t) - q^N(t)| \right) \to 0, \text{ as } N \to +\infty, \forall i, j \leq n.$$

Using the convergence *Lemma* 5 we then deduce that $u^N(t)$ converges uniformly to $u(t)$ in the norm of $L^2(\Omega)$, as $N \to +\infty$.

On the other hand, let $\bar{\mathcal{A}}^N(t)$ be the second order elliptic operator with the following constant average coefficients:

$$\bar{a}_{i,j}^N = \frac{1}{N} \sum_{k=0}^{N-1} a_{i,j}(t_k), \bar{b}_i^N = \frac{1}{N} \sum_{k=0}^{N-1} b_i(t_k) \text{ and } \bar{q}^N = \frac{1}{N} \sum_{k=0}^{N-1} q(t_k), \quad \forall i, j \leq n,$$

and let $\bar{u}^N$ be the solution of the corresponding parabolic equation with initial condition equal to $g$, i.e. solution of the parabolic equation (7) with the operator $\mathcal{A}(t)$ in the equation replaced by $\bar{\mathcal{A}}^N(t)$. We have that

$$|\bar{a}_{i,j} - \bar{a}_{i,j}^N| + |\bar{b}_i - \bar{b}_i^N| + |\bar{q} - \bar{q}^N| \to 0, \text{ as } N \to +\infty, \quad \forall i, j \leq n,$$

where the coefficients above are defined in (15). Applying the convergence *Lemma* 5 again we deduce that $\bar{u}^N(t)$ converges uniformly to $\bar{u}_T(t)$ in the norm of $L^2(\Omega)$, as $N \to +\infty$, where $\bar{u}(t)$ is solution of the parabolic problem in (16).

Now, using the result established in step 1 we have $\bar{u}^N(T) = u^N(T)$. Taking the limit, as $N \to +\infty$, we then deduce that $u(T) = \bar{u}_T(T)$, which completes the proof in step 2.

**Step 3**: we now prove the result when $g \in L^{2,s}(\Omega)$, for $s$ not necessarily equal to 0. From *Theorem* 1 there exists unique solutions $u$ and $\bar{u}_T$ that belong to $L^2\left(0,T;H_0^{1,s}(\Omega)\right)$ of the parabolic problems (7) and (16), respectively. We will show that $u(T) = \bar{u}_T(T)$. Let us denote by $C_0^\infty(\Omega)$ the set of infinitely differentiable function with compact support in the open set $\Omega$. Since $C_0^\infty(\Omega) \subset H_0^1(\Omega) \cap H^2(\Omega)$ is dense in $L^{2,s}(\Omega)$, there exist a sequence $g_N, N \geq 1$, that belongs to $H_0^1(\Omega) \cap H^2(\Omega)$ and converges to $g$ in the norm of $L^{2,s}(\Omega)$. Let $u_N$ and $\bar{u}_{T,N} \in L^2\left(0,T;H_0^{1,s}(\Omega)\right)$ be the respective solutions to the parabolic problems (7) and (16), where $g$ is replaced by $g_N$. From Step 2, and since $g_N \in H_0^1(\Omega) \cap H^2(\Omega)$, we have $u_N(T) = \bar{u}_{T,N}(T)$.

On the other hand, since $g_N$ converges to $g$ in $L^{2,s}(\Omega)$, as $N \to \infty$, we have that $u_N(T)$ converges to $u(T)$, and $\bar{u}_{T,N}(T)$ converges to $\bar{u}_T(T)$, in the norm of $L^{2,s}(\Omega)$ according to the energy estimates in *Theorem 1*. The proof in step 3, and that of *Theorem 2*, then completes by taking the limit in the equality $u_N(T) = \bar{u}_{T,N}(T)$, as $N \to +\infty$, which yields $u(T) = \bar{u}_T(T)$ as needed. **Q.E.D.**

## 5. Appendix: more Finance background of the B&S equation, and two different proofs when $\Omega = \mathbb{R}_+^n$

This appendix provides to the more interested reader more finance background on the B&S equation (1), and the knock out feature of financial options. It also presents for the sake of completeness two different proofs of the main result in the case where $\Omega' = \mathbb{R}_+^n$, which are not valid in the case where $\Omega' \neq \mathbb{R}_+^n$.

**Finance interpretation of the B&S equation (1):** The latter plays an important role in financial risk management, and is used in pricing and hedging financial derivative instruments, see [5] and [9] for more details. We assume the market to be complete and arbitrage free. The market is affected by $n$ sources of randomness, and there is available $n$ tradable assets with prices at time $\tau$ denoted by $S_i(\tau), i = 1, \ldots, n$. Also denote $S(\tau) = (S_1(\tau), \ldots, S_n(\tau))$. We assume the latter to be a multivariate lognormal stochastic differential process under a risk-neutral measure $\mathbb{Q}$:

$$dS_i(\tau) = (r(\tau) - m(\tau))S_i d\tau + \sigma_i(\tau)S_i dW_i, \quad i = 1, \ldots, n. \tag{31}$$

For each $i$, $\sigma_i(\tau)$ is the volatility parameter of $S_i(\tau)$, and $W_i(\tau), i = 1, \ldots, n$, are correlated Brownian motions under the risk-neutral measure $\mathbb{Q}$ (see [2], Chapter 6, on multi-asset models). The parameters $r$, $m$ are, respectively, the short term risk-free interest rate and the dividend yield. For each $i,j$, we denote by $\rho_{i,j}(\tau)$ the correlation parameter between $W_i(\tau)$ and $W_j(\tau)$. The usual shorthand notation for expressing this correlation relation is: $dW_i dW_j = \rho_{ij} d\tau, \forall i,j$, with $\rho_{i,i} = 1$.

Given the assumptions above, the price $U(\tau, y)$, at time $\tau \leq T$, and given that $y = S(\tau)$, of a derivative instrument that pays-off $g(S(T))$ at maturity $T$, satisfies the B&S equation (1), with $d = 0$, and with final condition $U(T, y) = g(y)$, see [5], ] and [9], where $g$ is a certain function defined on $\Omega$ characterizing the payoff of the derivative instrument. Some of the known derivative instruments are European put and call options. The payoff function corresponding to a European put option with strike price $K$ on the basket of assets $S_i(\tau), i = 1, \ldots, n$, is $g(y_1, \ldots, y_n) = \max(K - y_1 - y_2 - \cdots - y_n, 0)$.

More recently, Insurance companies have been selling Variable Annuity products, with embedded investment guarantees, such as Guaranteed Minimum Maturity Benefits (GMMB). The latter are put options in essence, with the difference that the contract is subject to policy decrements: essentially mortality and lapse decrements. If we denote the decrement rate to be $d$ then upon maturity $T$ of the contract the insurance company will be liable for the following contingent policy claim: $g(S_1(T), \ldots, S_n(T))$. In this case we have assumed that the fund manager follows a passive buy-and-hold investment strategy where the fund, with value $S_1(\tau) + \cdots + S_n(\tau)$, is not rebalanced. Assuming that the fund manager continuously rebalances the fund, so to keep constant proportions between the asset holdings $S_i(\tau), i = 1, \ldots, n$, would simplify the B&S equation (1) to a one space variable equation (i.e. $n = 1$), which is a much simpler problem and a particular case of the one under consideration.

Using the risk neutral valuation principle, the Fair Market Value $U(\tau, y)$ of the GMMB embedded in Variable Annuity contracts is given by the following risk neutral expectation:

$$U(\tau, y) = e^{-\int_\tau^T (r(s)+d(s))ds} \mathbb{E}^{\mathbb{Q}}_{S(\tau)=y} \left( g(S_1(T), \ldots, S_n(T)) \right).$$

It is then given by $U(\tau, y) = e^{-\int_\tau^T d(s)ds} U'(\tau, y)$, where $U'\tau, y)$ is known to satisfy the B&S Equation (1) with $d = 0$. Direct calculations then show that $U(\tau, y)$ satisfies equation (1) with the decrement factor $d$. The final condition $U(T, y) = g(y)$ also results. Equation (1) can also be directly established by applying the Feynman-Kac "with killing" theorem, see ], to the equation above. For GMMB the payoff at maturity is $g(y_1, \ldots, y_n) = \max(K - y_1 - \cdots - y_n, 0)$, where $K$ is the Guaranteed Minimum Maturity Benefits value. The parameter $m$ in the case of a variable annuity product is now interpreted differently and represents the Margin Expense Ratio (MER), which measures in practice the rate of daily deduction from the fund to cover for expenses, profit margins and other charges.

It is worth noting that the uniform elliptic property (2) is a consequence of the assumption that the market is complete and the coefficients $\sigma_i$ and $\rho_{i,j}$, $i, j \leq n$, which do not depend on $y$, are piecewise continuous with respect to $\tau$. In fact, for $S_i(\tau) = \xi_i, i \leq n$, for a given $\xi \in \mathbb{R}^n$, we (formally) have $\sum_{i,j \leq n} \rho_{i,j} \sigma_i \sigma_j \xi_i \xi_j \sqrt{d\tau} = $ VAR $(dS(\tau)) > 0$, $\forall \tau \geq 0$, $\forall \xi \in \mathbb{R}^n - \{0\}$, since the market is complete. Taking the minimum over the compact set of $(\tau, \xi) \in \mathbb{R}^{n+1}$, such that $\tau \in [0, T]$ and $|\xi| = 1$, we have that

$$\sum_{i,j \leq n} \rho_{i,j}(\tau) \sigma_i(\tau) \sigma_j(\tau) \xi_i \xi_j \geq c > 0, \quad \forall \xi \in \mathbb{R}^n \text{ with } |\xi| = 1, \text{ and } \forall \tau \in [0, T].$$

**On the Knock-out feature and relation to the Dirichlet condition:** The Knock-out and knock-in features are characteristics of a class of exotic options, called barrier options. A Knock-out feature causes the option to cease if the vector of asset values $S(\tau)$ reaches a certain boundary limit. An example of a knock-out option is an Up-and-Out option, which ceases if the total portfolio value $S_1(\tau) + \cdots + S_n(\tau)$ increases to or above a certain level $L$. The price of knock-out options still satisfy the B&S equation (1), but with the added condition of Dirichlet on the boundary that defines the knock-out feature. In the previous example, the Dirichlet condition would be imposed on the part of the boundary defined by $y_1 + y_2 + \cdots + y_n = L$ that is strictly included in $\mathbb{R}^n_+$. This is the part of the boundary $\partial \Omega'$ that do not intersect with $\partial \mathbb{R}^n_+$. The latter is in fact not reachable by any vector of multivariate lognormal variables, which is the model assumed in deriving the B&S equation (1) as explained above.

**Two direct proofs of the main result in the case of $\Omega' = \mathbb{R}^n_+$:** We have the following proposition:

***Proposition***: *We assume the coefficients $r, m, d$, $\sigma_i$ and $\rho_{i,j}$, $i, j \leq n$, to be piecewise continuous with respect to $\tau$ and independent of $y$. For $= \mathbb{R}^n_+$, and under the elliptic property (2), we have $\bar{U}(\tau_0, y) = U(\tau_0, y)$, for any given $\tau_0$, where $U$ and $\bar{U}$ are defined in in (1) and (3), respectively.*

We present here two straightforward proofs for the sake of completeness. The first proof is based on the Risk neutral valuation argument, usually used to solve problems arising in mathematical finance. The second proof is based on the Fourier transform, which cannot be used in the general case where $\Omega' \neq \mathbb{R}^n_+$.

**Proof 1 (Risk Neutral valuation approach):** From (31) the stochastic variables $\ln(S_i)$ satisfy $d\ln(S_i) = \left(r(\tau) - m(\tau) - \frac{\sigma_i^2(\tau)}{2}\right)dt + \sigma_i(\tau)dW_i$, which after integrating between $\tau_0$ and $T$ shows that $(\ln(S_1(T)), \ldots, \ln(S_n(T)))$ is a multivariate normal distribution with a mean equal to the vector with components $\ln(S_i(\tau_0)) + \left(\bar{r} - \bar{m} - \frac{\bar{\sigma}_i^2}{2}\right)(T - \tau_0), i = 1, \ldots, n$, and variance-covariance matrix with components equal to $\int_{\tau_0}^T \sigma_i \sigma_j \rho_{i,j} ds = \bar{\rho}_{i,j} \bar{\sigma}_i \bar{\sigma}_j (T - \tau_0), i, j \leq n$. Note that the variance-covariance matrix, which is symmetric, is also invertible for each $\tau$ since the market is assumed to be complete. The parameters $\bar{r}, \bar{m}, \bar{\sigma}_i^2$ and $\bar{d}$ are average values as defined in (4).

Consider now $\overline{W}_i(\tau), i = 1, \ldots, n$, to be correlated Brownian motions with constant correlations $\bar{\rho}_{i,j}$ equal to (5), and the variables $\bar{S}_i(\tau)$ solutions to the following stochastic differential equations:

$$d\bar{S}_i(\tau) = (\bar{r} - \bar{m})\bar{S}_i d\tau + \bar{\sigma}_i \bar{S}_i d\overline{W}_i, \quad \text{with } \bar{S}_i(\tau_0) = S_i(\tau_0), \text{ for } i = 1, \ldots, n.$$

As for (31) we have that $(\ln(\bar{S}_1(T)), \ldots, \ln(\bar{S}_n(T)))$ is also a correlated multivariate normal distribution with the same mean and variance-covariance matrix as that of $(\ln(S_1(T)), \ldots, \ln(S_n(T)))$. This means that the two vector variables have identical probability distributions, which implies that they have the same expectation value:

$$U(\tau_0, y) = e^{-\int_{\tau_0}^T (r(s) + d(s))ds} \mathbb{E}^{\mathbb{Q}}_{S(\tau_0)=y}\left(g(S_1(T), \ldots, S_n(T))\right) = e^{-(\bar{r}+\bar{d})(T-\tau_0)} \mathbb{E}^{\mathbb{Q}}_{\bar{S}(\tau_0)=y}\left(g(\bar{S}_1(T), \ldots, \bar{S}_n(T))\right).$$

Let us denote by $V(\tau, y)$ the function defined as in the last term in the equation above:

$$V(\tau, y) \stackrel{\text{def}}{=} e^{-(\bar{r}+\bar{d})(T-\tau)} \mathbb{E}^{\mathbb{Q}}_{\bar{S}(\tau)=y}(g(\bar{S}_1(T), \ldots, \bar{S}_n(T))).$$

According to the Feynman-Kac theorem, with killing, this function above satisfies equation (3) with the final condition $V(T, y_1, \ldots, y_n) = g(y_1, \ldots, y_n)$, which proves that $U(\tau_0, y) = V(\tau_0, y) = \overline{U}(\tau_0, y)$, assuming existence and uniqueness of solutions of the equation (3).

**Proof 2 (using the Fourier transform):** Consider the following Fourier transform with respect to $x$ of $u(t, x) = U(\tau, y)$, where $x_i = \text{Ln}(y_i), i = 1, \ldots, n$, and $t = (T - \tau)$:

$$v(t, \xi) = \int_{\mathbb{R}^n} u(t, x) e^{-2\omega\pi\xi x} dx, \quad \forall \xi \in \mathbb{R}^n, \text{ with } \omega^2 = -1.$$

Equation (6) then yields $\frac{dv}{dt} + P(t, \xi)V = 0$, which is a linear ordinary differential equation of the first order for each $\xi \in \mathbb{R}^n$. Here $P(t, \xi)$ is a second order polynomial with respect to $\xi$ equal to:

$$P(t, \xi) = \frac{(2\pi)^2}{2} \sum_{i,j=1}^n \rho_{i,j}\sigma_i\sigma_j\xi_i\xi_j + 2\omega\pi \sum_{i=1}^n \left(\frac{\sigma_i^2}{2} - (r - m)\right)\xi_i + r + d. \tag{32}$$

The solution of this differential equation, for each $\xi \in \mathbb{R}^n$, is given by $v(t, \xi) = e^{-\int_0^t P(s,\xi)ds} v(0, \xi)$, so that $v(T - \tau_0, \xi) = \bar{v}(T - \tau_0, \xi) = e^{-(T-\tau_0)\bar{P}(\xi)} v(0, \xi)$, where

$$\bar{P}(\xi) = \frac{(2\pi)^2}{2} \sum_{i,j=1}^n \bar{\rho}_{i,j}\bar{\sigma}_i\bar{\sigma}_j\xi_i\xi_j + 2\omega\pi \sum_{i=1}^n \left(\frac{\bar{\sigma}_i^2}{2} - (\bar{r} - \bar{m})\right)\xi_i + \bar{r} + \bar{d}.$$

Here, the coefficients of the polynomial $\bar{P}(\xi)$ are defined in (4) and (5). The proof concludes by noting that $\xi \to \bar{v}(T - \tau_0, \xi)$ is the Fourier transform, with respect to $x$, at $t = T - \tau_0$ of the solution of the equivalent second order parabolic problem where the coefficients are given by the constants defined in (4) and (5). **Q.E.D.**